\documentclass[10pt]{article}
\usepackage[latin1]{inputenc}
\usepackage[english]{babel}
\usepackage{amssymb,amsfonts, amsmath}
\usepackage{amsthm}

\usepackage{a4wide}
\newtheorem{theorem}{Theorem}[section]
\newtheorem{proposition}[theorem]{Proposition}

\newtheorem{lemma}[theorem]{Lemma}
\newtheorem{remark}[theorem]{Remark}
\newtheorem{definition}[theorem]{Definition}
\newtheorem{example}[theorem]{Example}

\newcommand{\bcl}{\begin{center}}
\newcommand{\ecl}{\end{center}}
\newcommand{\brl}{\begin{right}}
\newcommand{\erl}{\end{right}}
\newcommand{\ben}{\begin{enumerate}}
\newcommand{\een}{\end{enumerate}}
\newcommand{\barr}{\begin{array}}
\newcommand{\earr}{\end{array}}
\newcommand{\btab}{\begin{tabular}}
\newcommand{\etab}{\end{tabular}}
\newcommand{\bdoc}{\begin{document}}
\newcommand{\edoc}{\end{document}}
\newcommand{\beqy}{\begin{eqnarray}}
\newcommand{\eeqy}{\end{eqnarray}}

\newcommand{\beqi}{\begin{eqnarray*}}
\newcommand{\eeqi}{\end{eqnarray*}}
\newcommand{\bitem}{\begin{itemize}}

\newcommand{\eitem}{\end{itemize}}
\newcommand{\nln}{\newline}
\newcommand{\newt}{\newtheorem}

\newcommand{\pa}{\partial}
\newcommand{\re}{{I\!\!R}}
\newcommand{\ren}{\re^N}
\newcommand{\xr}{x\in\re }
\newcommand{\x}{\times}
\newcommand{\dyle}{\displaystyle}
\newcommand{\ene}{{I\!\!N}}
\newcommand{\irn}{\int\limits_{\re^N}}
\newcommand{\io}{\int\limits_{\O}}
\newcommand{\meas}{{\rm meas\,}}

\newcommand{\sign}{{\rm sign}}
\newcommand{\map}{\longrightarrow }
\newcommand{\imp}{\Longrightarrow }
\renewcommand{\div}{\nabla\cdot }
\newcommand{\sen}{{\rm sen\,}}
\newcommand{\tg}{{\rm tg\,}}
\newcommand{\arcsen}{{\rm arcsen\,}}
\newcommand{\arctg}{{\rm arctg\,}}
\newcommand{\supp}{{\textsl supp\ }}
\newcommand{\ity}{\int_{-\iy}^{+\iy}}
\newcommand{\limit}{\lim\limits}
\newcommand{\limi}{\limit_{n\to\infty}}
\newcommand{\sumi}{\sum\limits_{n=1}^{\infty}}
\newcommand{\ulu}{\underline u}
\newcommand{\ulw}{\underline w}
\newcommand{\ulz}{\underline z}
\newcommand{\ulv}{\underline v}
\newcommand{\uls}{\underline s}
\newcommand{\olu}{\overline u}
\newcommand{\olv}{\overline v}
\newcommand{\ols}{\overline s}
\newcommand{\ob}{\overline\b}
\newcommand{\ovar}{\overline\var}
\newcommand{\wv}{\widetilde v}
\newcommand{\wu}{\widetilde u}
\newcommand{\ws}{\widetilde s}
\renewcommand{\a }{\alpha }
\renewcommand{\b }{\beta }
\newcommand{\g }{\gamma}
\newcommand{\G }{\Gamma }
\renewcommand{\d }{\delta }

\newcommand{\D }{\Delta }
\newcommand{\e }{\varepsilon }
\newcommand{\z }{\zeta }
\renewcommand{\l }{\lambda }
\renewcommand{\L }{\Lambda }
\newcommand{\m }{\mu }
\newcommand{\n }{\nabla }
\newcommand{\s }{\sigma }
\newcommand{\Sig }{\Sigma }
\renewcommand{\t }{\tau }
\newcommand{\var }{\varphi }
\renewcommand{\o }{\omega }
\renewcommand{\O }{\Omega }
\newcommand{\bR}{{\bf R}}
\newcommand{\bC}{{\bf C}}
\newcommand{\bZ}{{\bf Z}}
\newcommand{\bN}{{\bf N}}
\newcommand{\bQ}{{\bf Q}}
\newcommand{\bK}{{\bf K}}
\newcommand{\bI}{{\bf I}}
\newcommand{\bv}{{\bf v}}
\newcommand{\bV}{{\bf V}}

\newcommand{\pair}[1]{g\left(#1\right)}

\def\qed{\unskip\kern 6pt \penalty 500
\raise -2pt\hbox{\vrule \vbox to10pt{\hrule width 4pt
\vfill\hrule}\vrule}\par}

\newenvironment{Proof}{\removelastskip\vskip12pt
plus 1pt \noindent\em\rm}{\hfill {\qed \hskip .2cm}}

\title
{Global Solutions of Semilinear Parabolic Equations on Negatively Curved Riemannian Manifolds}
\author{Fabio Punzo\thanks{Dipartimento di Matematica, Politecnico di Milano, Italia (fabio.punzo@polimi.it).}}
\date{}

\begin{document}
 \maketitle
{\abstract{ \noindent We are concerned with global existence for semilinear parabolic equations on Riemannian manifolds with negative sectional curvatures. A particular attention is paid to the class of initial conditions which ensure existence of global solutions. Indeed, we show that such a class is crucially related to the curvature bounds.
\bigskip


}}

\bigskip
\medskip
\smallskip

\section{Introduction} \setcounter{equation}{0}
We investigate existence of global nonnegative
solutions of Cauchy problems for semilinear parabolic equations of the following
form:
\begin{equation}\label{e1}
\left\{
\begin{array}{ll}
 \,  \pa_t u = \Delta u\, +\, h(t) u^p \, &\textrm{in}\,\,M\times (0,T)
\\&\\
\textrm{ }u \, = u_0& \textrm{in}\,\,  M\times \{0\} \,;
\end{array}
\right.
\end{equation}
here $M$ is a complete noncompact
Riemannian manifold of dimension $n$, endowed with a metric tensor $g$, with strictly negative sectional
curvatures, $\Delta$ is the Laplace-Beltrami operator on $M$; $h$
is a positive continuous function defined in $[0,\infty)$, the
initial datum $u_0$ is continuous, nonnegative and bounded on $M$, $p>1$\,.

Problem \eqref{e1} with $M=\mathbb R^n$ and $h(t)\equiv 1$ has been the object of detailed investigations. In particular, it is well-known that for any $u_0\not\equiv 0$, the solution blows-up in finite time, if $1<p\le 1+\frac 2
n\,$ (see \cite{Fuji}, \cite{Lev}). On the contrary, if $p>1+\frac 2 n$ and the initial datum $u_0$ is small enough, then
global solutions exist.

Problem \eqref{e1} with $M$ the hyperbolic space $\mathbb H^n$ has been addressed in \cite{BPT}. It is shown that
if $h(t)\equiv 1\; (t\ge 0)$ or
\begin{equation}\label{e56}
\a_1 t^q\le h(t) \le \a_2 t^q\quad \big(t\geq t_0\big)
\end{equation}
for some $\a_1>0, \a_2>0, t_0>0$ and $q>-1\,$,
then {\it for any} $p>1$ there
exist global solutions for sufficiently small initial data $u_0$. Moreover, if
$h(t)=e^{\s t}\; (\s >0)$, we can have both finite time blow-up and global existence, in dependence of the parameter $\s$. In fact, let $\l_1(\mathbb H^n)$ be the infimum of the $L^2-$ spectrum of the operator
$-\Delta\,$ on $\mathbb H^n$\,. If $1< p< 1+\frac{\s}{\l_1(\mathbb H^n)},$ then every nontrivial solution of problem \eqref{e1} blows-up in finite time; while (see also \cite{WY} if $p \geq 1+\frac{\s}{\l_1(\mathbb H^n)},$ then problem \eqref{e1} has global solutions for small initial data.

\medskip

The results given in \cite{Fuji} have been extended to general Riemannian
manifolds $M$ in \cite{Zhang} and in \cite{MMP}, under suitable growth conditions on volume of geodesic balls. However, such hypotheses cannot be satisfied on Riemannian manifolds with strictly negative sectional curvatures. On the other hand, in \cite{Pu12} the results established in \cite{BPT} have been generalized to Cartan-Hadamard manifolds $M$ with sectional curvatures
bounded above by a negative constant. Some global existence results for mild solutions belonging to $C([0, T); L^p(M))$ have been established in \cite{Pu14}, by using general results in semigroup theory stated in \cite{W}. We should point out that in the present paper we shall use complete different methods. Indeed, we consider  bounded weak solutions; moreover, our assumptions on initial datum will be completely different in character and related to the geometry of the underlying manifold, and that on the function $h$ are weaker.  

Both in \cite{BPT} and in \cite{Pu12}  the elliptic equation
\begin{equation}\label{e22}
\Delta \phi + \l \phi\, =\, 0\quad \textrm{in}\;\; M\,
\end{equation}
has a prominent role. In fact, it is well-known that for any $\l\le \l_1$ there exists a
classical positive solution $\phi$ to equation \eqref{e22} (see
\cite{CYau}, \cite{Grig4}). When $\l=\l_1,$ we say that $\phi$ is a {\it ground state} on $M$\,. However, in general, it is not known the behaviour of $\phi$ at infinity and whether $\phi$ belongs to some Lebesgue space. Indeed, in \cite{BPT} it is shown that $\phi\in L^\infty(\mathbb H^n)\setminus L^2(\mathbb H^n).$ Instead, in \cite{Pu12} it is directly assumed that $\phi\in L^{\infty}(M)$ for $\l=\l_1.$ By means of a bounded solution of equation \eqref{e22} with $\l=\l_1$, both in \cite{BPT} and in \cite{Pu12}, global existence for problem \eqref{e1} is deduced. We should mention that in \cite{Donn} and in \cite{CiaMa} sufficient conditions for the boundedness of $\phi$ are given. Observe that in \cite{Donn} it is supposed that $\phi\in L^2(M)$, while in \cite{CiaMa} the case $\mu(M)<\infty$ is addressed. However, in our situation $\mu(M)=\infty$. Moreover, we do not know in general whether $\phi\in L^2(M)$; indeed, as recalled above, in the special case $M=\mathbb H^n$ it is not true. Hence, we cannot use the results stated in \cite{CiaMa}, \cite{Donn}\,.
\smallskip

In this paper we always assume that the sectional curvatures are bounded above by a negative constant; moreover, we suppose that the radial sectional curvature negatively diverges at infinity with a certain rate (see assumption \eqref{A0} below). However, for some results, the latter assumption is not required. Such a class of Cartan-Hadamard manifolds have been recently considered also in \cite{GMP1}, \cite{GMP}, \cite{GMV}, for different purposes. Our goal is to study global existence for problem \eqref{e1}, without assuming that a ground state is bounded. Furthermore, we aim at showing that the class of initial conditions for which global existence prevails is related to the bounds on the sectional curvatures.

We consider bounded weak supersolutions of equation \eqref{e22}, for each $\l\in (0, \l_1]$. More precisely, we prove that for any $\l\in (0, \l_1]$  a bounded weak supersolution $w$ of equation \eqref{e22} exists; in addition, $w$ decays at infinity. Using such a supersolution we can prove that problem \eqref{e1} admits global solutions for any $p>1$, when $h(t)\equiv 1$ or \eqref{e56} holds, provided that $u_0\leq C w$ for a properly chosen positive constant $C$. Moreover, let $\l_1(M)$ be the infimum of the $L^2-$ spectrum of the operator
$-\Delta\,$ on $M$. We obtain global existence also for $h(t)=e^{\s t}\, (\s>0)$ and $p > 1+\frac{\s}{\l_1(M)}.$ On the other hand, for $p < 1+\frac{\s}{\l_1(M)},$ there is finite time blow-up, by results in \cite{Pu12}.

Observe that, in order to deal with special functions $h$, we prove that there exist $\l\in (0, \l_1]$ and a weak bounded supersolution of  \eqref{e22} for such $\l$. Then we obtain global existence for problem \eqref{e1}. Hence it is not necessary to take $\l=\l_1.$ Note that a special choice of $h$ is for instance $h\equiv 1$; so, we can consider the model problem
\begin{equation}\label{e1b}
\left\{
\begin{array}{ll}
 \,  \pa_t u = \Delta u\, +\, u^p \, &\textrm{in}\,\,M\times (0,T)
\\&\\
\textrm{ }u \, = u_0& \textrm{in}\,\,  M\times \{0\} \,.
\end{array}
\right.
\end{equation}
Moreover, for small values of $\l<\l_1$, we can construct bounded supersolutions of equation \eqref{e22} which are more explicit than $w$ introduced above. Consequently, we can make more explicit assumptions on initial conditions $u_0$ ensuring global existence for problem \eqref{e1b}. More precisely, we see that such supersolutions decay at infinity with a rate which depends on the behaviour at infinity of the sectional radial curvature. In particular, if the sectional radial curvature negatively diverges fastly at infinity, then the supersolutions decay slowly, and thus we can enlarge the class of initial data $u_0$ in \eqref{e1b}.
Concerning the function $h$, besides $h\equiv 1$, for instance, we can also suppose that $h$ fulfills \eqref{e56}. However, if we take $h(t)=e^{\s t}\, (\s>0)$, we have to require a stronger hypothesis on $p$ than the previous one associated with the supersolution $w$ considered above. To be specific, we must assume that $p > 1+\frac{\s}{\l},$ and obviously, for small values of $\l<\l_1$, $1+\frac{\s}{\l}>1+\frac{\s}{\l_1(M)}.$

\smallskip

The paper is organized as follows. In Section \ref{MB} we recall some useful preliminaries from Riemannian Geometry and we introduce our geometric assumptions. In Section \ref{MR} we state our main results and we discuss some examples. Section \ref{supersolutions} is devoted to the construction of various supersolutions of equation \eqref{e22}; this can also have a certain independent interest. Finally, in Section \ref{dim1} we give the proofs of the global existence results.

\section{Mathematical framework}\label{MB} \setcounter{equation}{0}
\subsection{Preliminaries from Riemannian Geometry}
In this section we collect some useful notions and results from Riemannian Geometry (see e.g. \cite{Grig}).

Let $M$ be a complete noncompact Riemannian manifold of dimension $n$. Let $\Delta$
denote the Laplace-Beltrami operator, $\nabla$ the Riemannian
gradient and $d\mu$ the Riemannian volume element on $M$.
In the sequel we consider {\it Cartan-Hadamard} manifolds, i.e. simply connected complete noncompact Riemannian manifolds with
nonpositive sectional curvatures. Observe that (see, e.g. \cite{Grig}) for a Cartan-Hadamard manifold $M$
the {\it cut locus } of $o$, $\operatorname{Cut}(o)$, is empty for any point $o \in M,$ thus $M$ is a manifold with a pole.
For any $x\in M\setminus \{o\}$, one can define the {\it polar
coordinates} with respect to $o$. Namely, for
any point $x\in M\setminus\{o\}$
there correspond a polar radius $r(x) :=\operatorname{dist}(x, o)$
and a polar angle $\theta\in \mathbb S^{n-1}$ such that the shortest
geodesics from $o$ to $x$ starts at $o$ with direction $\theta$ in
the tangent space $T_o M$. Since we can identify $T_o M$ with
$\mathbb R^n$, $\theta$ can be regarded as a point of $\mathbb
S^{n-1}.$ For any $x_0\in M$ and for any $R>0$ we set $B_R(x_0):=
\big\{x\in M\,:\, \operatorname{dist}(x, x_0)<R\,\big\}.$

The Riemannian metric in $M\setminus\{o\}$ in polar coordinates reads
\[g= dr^2+A_{ij}(r, \theta)d\theta^i d\theta^j, \]
where $(\theta^1, \ldots, \theta^{n-1})$ are coordinates in $\mathbb S^{n-1}$ and $(A_{ij})$ is a positive definite matrix.
It is not difficult to see that the Laplace-Beltrami operator in polar coordinates has the form
\begin{equation}\label{e6} \Delta = \frac{\partial^2}{\partial r^2} +
\mathcal F(r, \theta)\frac{\partial}{\partial r}+\Delta_{S_{r}},
\end{equation}
where $\mathcal F(r, \theta):=\frac{\partial}{\partial
r}\big(\log\sqrt{A(r,\theta)}\big)$, $A(r,\theta):=\det
(A_{ij}(r,\theta))$, $\Delta_{S_r}$ is the Laplace-Beltrami operator
on the submanifold $S_{r}:=\partial B_r(o)$\,.

A manifold with a pole is a {\it spherically symmetric manifold} or a {\it model}, if the Riemannian metric is given by
\begin{equation}\label{e7}
g= dr^2+\psi^2(r)d\theta^2,
\end{equation}
where $d\theta^2=\beta_{ij}d\theta^i d \theta^j$ is the standard
metric in $\mathbb S^{n-1}$, $\beta_{ij}$ being smooth functions of
$\theta^1, \ldots, \theta^{n-1},$ and $\psi\in \mathcal A$, where
$$\mathcal A:=\left\{f\in C^\infty((0,\infty))\cap C^1([0,\infty)): \, f'(0)=1, \, f(0)=0, \, f>0 \ \textrm{in}\;\, (0,\infty)\right\} .$$
In this case, we write $M\equiv
M_\psi$; furthermore, we have $\sqrt{A(r,\theta)}=\psi^{n-1}(r)$, so
that
\[\Delta = \frac{\partial^2}{\partial r^2}+ (n-1)\frac{\psi'}{\psi}\frac{\partial}{\partial r}+ \frac1{\psi^2}\Delta_{\mathbb S^{n-1}}\,,\]
where $\Delta_{\mathbb S^{n-1}}$ is the Laplace-Beltrami operator in
$\mathbb S^{n-1}\,.$ Observe that for $\psi(r)=r$, $M=\mathbb R^n$, while for $\psi(r)=\sinh r$, $M$ is the $n-$dimensional hyperbolic space $\mathbb H^n$.

\smallskip

Let us recall comparison results for sectional and Ricci curvatures
that will be used in the sequel. For
any $x\in M\setminus\{o\}$, denote by
$\textrm{Ric}_o(x)$ the {\it Ricci curvature} at $x$ in the
direction $\frac{\partial}{\partial r}$. For every $x\in M$ and for every plane $\pi\subseteq T_x M$ denote
by $K_{\pi}(x)$ the {\it sectional curvature} of the plane $\pi$
(see \cite{GHL}). Let $\omega$ denote any
pair of tangent vectors from $T_xM$ having the form
$\left(\frac{\partial}{\partial r} ,X\right)$, where $X$ is a unit
vector orthogonal to $\frac{\partial}{\partial r}$. Denote by
$\textrm{K}_{\omega}(x)$ the {\it sectional curvature} at the point
$x$ of the $2$-section determined by $\omega$; it is also called {\it sectional radial curvature}. Observe that (see
\cite[Section 15]{Grig}), if
\begin{equation}\label{e9}
\textrm{K}_{\omega}(x)\leq -\frac{\psi''(r)}{\psi(r)}\quad \textrm{for all}\;\; x=(r,\theta)\in M\setminus\{o\},
\end{equation}
for some function $\psi\in \mathcal A$, then
\begin{equation}\label{e10}
\mathcal F(r, \theta)\geq (n-1)\frac{\psi'(r)}{\psi(r)}\quad
\textrm{for all}\;\; r>0,\, \theta \in \mathbb S^{n-1}\,.
\end{equation}

On the other hand, if
\begin{equation}\label{e11}
\textrm{Ric}_{o}(x)\geq -(n-1)\frac{\phi''(r)}{\phi(r)}\quad \textrm{for all}\;\; x=(r,\theta)\in M\setminus\{o\},
\end{equation} for some function $\phi\in \mathcal A$, then
\begin{equation}\label{e12}
\mathcal F(r, \theta)\leq (n-1)\frac{\phi'(r)}{ \phi(r)}\quad \textrm{for all}\;\; r>0, \theta \in \mathbb S^{n-1}\,\,\textrm{with}\,\, x=(r,
\theta)\in M\setminus \{o\}\,.
\end{equation}
If $M_\psi$ is a model manifold, then for any $x=(r,
\theta)\in M_\psi\setminus\{o\}$
\[\textrm{K}_{\omega}(x)=-\frac{\psi''(r)}{\psi(r)},\]
and
\[\textrm{Ric}_{o}(x)=-(n-1)\frac{\psi''(r)}{\psi(r)}\,.\]

\subsection{Geometric assumptions and consequences}

Throughout the paper we shall make the following hypothesis:
\begin{equation} \tag{{\it $A_0$}} \label{A0}
\begin{cases}
\textrm{(i)} & M \textrm{ is a Cartan-Hadamard manifold of dimension $n\ge2$} \, ; \\
\textrm{(ii)} & \textrm{there exists}\;k>0\;\textrm{such that for any}\,x\in M \,\textrm{and for any}\\ & \textrm{plane}\;\pi\subseteq T_x M\,\,\textrm{there
holds}\;\, K_{\pi}(x)\le - k^2\,;\\
\textrm{(iii)} & \textrm{K}_\omega(x)\leq -C_0 \left(1+d(x,o)^\gamma \right) \ \textrm{for some } o\in M, C_0>0 \ \textrm{and } \gamma\geq 0 \, .
\end{cases}
\end{equation}
Note that, in particular, hypothesis \eqref{A0}-(ii) implies that
\begin{equation}\label{eq35}
K_\omega(x) \leq - k^2 \quad  \textrm{for any } x \equiv (r,\theta) \in M \setminus \{ o \} \, .
\end{equation}
For instance, assumption \eqref{A0} is satisfied if $M=\mathbb H^N$, with $ \gamma=0 $. More generally, it is not difficult to show that \eqref{A0} is fulfilled e.g. by Riemannian models associated with suitable convex functions $\psi\in \mathcal A$ such that
\begin{equation}\label{e50z}
\psi(r) = e^{f(r)} \, , \quad f(r) \sim C \, r^{1+\frac{\gamma}2}  \quad \textrm{as } r \to \infty \quad \textrm{if } \gamma>0\,,
\end{equation}
where $C$ are positive constants. Here by $ f(r) \sim g(r) $ we mean that the ratio $ f(r)/g(r) $ tends to $1$ as $ r \to \infty $.
We refer the reader to \cite[Section 2.3]{GMV} for more details in this regard.

An important role in the following will be played by the next lemma.
\begin{lemma}\label{lemma1}
Let assumption \eqref{A0} be satisfied. Then there exists a positive constant $ \underline C$, depending on $ C_0,\gamma$, such that
\begin{equation}\label{eq15}
\mathcal F(r,\theta) \ge \frac{\underline C(n-1)}{r}(1+r)^{1+\frac{\g}{2}} \quad \textrm{for any }\,  x \equiv (r,\theta) \in M \setminus \{ o \} \, .
\end{equation}
Furthermore,
\begin{equation}\label{eq15b}
\mathcal F(r, \theta)\geq k \quad \textrm{for any} \;\,  x \equiv (r,\theta) \in M \setminus \{ o \} \, .
\end{equation}
\end{lemma}
\begin{proof} Inequality \eqref{eq15} is shown in \cite{GMV} (see also \cite{GMP}), by using \eqref{e10} for a suitable choiche of $\psi\in \mathcal A$. Furthermore, \eqref{eq15b} follows from \eqref{e10}, by choosing $\psi(r)=\frac 1 k\sinh(k r)$ with $k$ given in \eqref{A0}-(ii) .
\end{proof}

\begin{remark}\label{oss2}
Note that in \cite{GMV} it is shown that, under assumption \eqref{A0}-(ii),
\[\mathcal F (r, \theta) \geq \hat C (n-1)r^{\frac{\g}{2}} \quad \textrm{for every}\;\; x\equiv (r, \theta)\in M\setminus B_{\hat R},\]
for some $\hat R>0$. Moreover, $\hat C\to \sqrt{C_0}\left(1 +\frac{\g}2\right)$ as $\hat R\to \infty$.
However, we do not know the precise value of $\underline C$ appearing in inequality \eqref{eq15}.
\end{remark}

\medskip

Let $spec(-\Delta)$ be the spectrum in $L^2(M)$ of the operator
$-\Delta$. Note that (see \cite{Grig4}, Chapter $4$)
\[spec(-\Delta)\subseteq [0,\infty)\,.\]
Denote by $\l_1(M)$ the bottom of $spec(-\Delta)$, that is
\[\l_1(M):=\inf spec(-\Delta)\,.\]

Observe that if \eqref{A0}(i)-(ii) holds, then (see \cite{McK}; see also
\cite{Grig})
\begin{equation}\label{e7}
\l_1(M)\ge \frac{(n-1)^2}4 k^2.
\end{equation}

\subsection{Definition of solutions}
We always make the following assumption:
\begin{equation} \tag{{\it $A_1$}} \label{A1}
\begin{cases}
\textrm{(i)} & h\in C([0,\infty))\,,\,h>0\,\, \textrm{in}\,[0,\infty)\, ; \\
\textrm{(ii)} & u_0\,\,\hbox{is continuous, nonnegative and bounded in}\;M\,.
\end{cases}
\end{equation}

We shall deal both with classical and with weak solutions to problem \eqref{e1} and to equation \eqref{e22}.
Weak solutions are meant in the following sense.
\begin{definition}\label{defwsol}
A {\em weak supersolution} to problem \eqref{e1} is a function $u\in
C\big(M\times [0, \t] \big)\cap L^\infty(M\times (0,\t))$ for any
$\t\in [0, T)$ such that
\begin{equation}\label{eq60}
\begin{aligned}
& -\int_0^\t \int_M u(x,t)\big\{\Delta\psi(x,t) +\pa_t \psi(x,t) \big\}d\mu dt \\ & \geq  \int_M u_0(x) \psi(x,0)d\mu +
\int_0^\t \int_M h(t)u^p(x) \psi(x,t) d\mu dt
\end{aligned}
\end{equation}
for any $\t\in [0, T)$, for any $\psi\in C^{2,1}(M\times
[0,\t]), \psi\geq 0$ with $supp\; \psi (\cdot,t)\subseteq M\;(t\in [0,\t])$
and $\psi(\cdot, \t)=0$\,.

{\em Weak subsolutions} to problem \eqref{e1} are defined similarly, replacing $\geq$ by $\leq$ in \eqref{eq60}. A {\em weak solution} is both a weak subsolution and a weak supersolution.
\end{definition}

\smallskip

\begin{definition}
A solution to problem \eqref{e1} is called {\em global}, if it
exists for any $t>0$, that is if \, $T=\infty$.

Instead, we say that a solution to problem \eqref{e1} {\em
blows-up in finite time}, if
\[\lim_{t\to T^-}\|u(\cdot, t)\|_{L^\infty(M)}\,=\,\infty\,,\]
for some $T>0\,.$
\end{definition}

\begin{definition}\label{solell}
A {\em weak supersolution} of equation \eqref{e22} is a function $u\in W^{1,2}_{\textrm{loc}}(M)$ such that
\begin{equation}\label{eq61} - \int_M \langle \nabla \psi(x), \nabla u(x) \rangle d\mu \leq \l \int_M u(x) \psi(x) d\mu
\end{equation}
for any $\psi\in C^1_{\textrm{c}}(M), \psi\geq 0.$

{\em Weak subsolutions} to problem \eqref{e22} are defined similarly, replacing $\geq$ by $\leq$ in \eqref{eq61}. A {\em weak solution} is both a weak subsolution and a weak supersolution.
\end{definition}

\section{Main Results}\label{MR} \setcounter{equation}{0}
Set $\l_1\equiv \l_1(M)$. For every $\l\in (0, \l_1]$ let
\[\tilde h(t):=h(t)e^{-(p-1)\l t}\quad \textrm{for any}\,\; t\ge 0,\]
\[\tilde H(t):= \int_0^t \tilde h(s) ds\quad \textrm{for any}\,\, t\ge 0\,,\]
Suppose that
\begin{equation}\label{e25}
\tilde H_{\infty}:=\lim_{t\to \infty}\tilde H(t)<\infty\,.
\end{equation}

For each $\alpha>0, \beta>0$, define
\begin{equation}\label{eq16}
v(x):= e^{-\b[r(x)]^{\alpha}}, \quad x\in M\,.
\end{equation}
We shall write $v(x)\equiv v(r).$

\smallskip

Under hypothesis \eqref{A0} with $\g>0$, for every $\l\in (0, \l_1]$ we shall construct a positive bounded weak supersolution $w$ of equation \eqref{e22} such that
\begin{equation}\label{eq1}
0<w(x)\leq  v(x) \quad \textrm{for all} \;\; x\in M\setminus B_{R_0}\,,
\end{equation}
for any $\b>0$, for any $\a>0$ satisfying
\begin{equation}\label{eq30}
\max\left\{1-\frac{\g}2, 0\right\}<\alpha<1+\frac{\g}2\,,
\end{equation}
for  $R_0>0$ sufficiently large (see Proposition \ref{prop1} below).

By using $w$, we can prove the following theorem.

\begin{theorem}\label{thm1}
Let assumptions \eqref{A0}-\eqref{A1} be satisfied with $\g>0$. Let $\l\in (0, \l_1].$
Suppose that condition \eqref{e25} is satisfied.
Moreover, assume that
\begin{equation}\label{eq3}
0\le u_0\le \tilde C w \quad \textrm{in} \;\; M\,,
\end{equation}
where
\begin{equation}\label{eq2}
0< \tilde C< \frac 1{\|w \|_{\infty}} \left[ \frac 1{(p-1)\tilde
H_{\infty}}\right]^{\frac 1{p-1}}\,.
\end{equation}
Then there exists a global solution $u$ of problem
\eqref{e1}; in addition, there exists $\bar C>0$ such
that
\begin{equation}\label{eq4}
\|u(\cdot, t)\|_{L^\infty(M)} \le \bar C\quad \textrm{for all}\;\;
t>0\,.
\end{equation}
\end{theorem}

If we only consider sufficiently small $\l>0$, we can make more explicit assumptions on the initial conditions $u_0$. This requires to know a bound from below on the supersolution of equation \eqref{e22}. In particular, we address the cases in which \eqref{A0} hold with $\g\geq 0 $ or $\g>0$ or $\g>2.$ We shall see that as $\g$ becomes bigger, we can enlarge the class of initial data $u_0$.

\smallskip

We show that, for every
\begin{equation}\label{eq32}
0<\l\leq \frac{(n-1)^2}{4} k^2\,,
\end{equation}
with $k>0$ given by \eqref{A0}-(ii), the function $v$ with

\begin{equation}\label{eq46}
\a=1\,,
\end{equation}
\and
\begin{equation}\label{eq46b}
\frac{k(n-1)-\sqrt{[k(n-1)]^2-4\l}}{2} \leq \b\leq \frac{k(n-1)+\sqrt{[k(n-1)]^2-4 \l}}{2}\,,
\end{equation}
is a weak supersolution of equation \eqref{e22} (see Proposition \ref{prop2} below). So, using such a $v$ we can prove the next result.

\begin{theorem}\label{thm2}
Let assumptions \eqref{A0}-\eqref{A1} be satisfied. Suppose that
conditions \eqref{eq32}, \eqref{e25} are satisfied, and that \eqref{eq3} and \eqref{eq2} hold with $w$ replaced by $v$.
Then there exists a global solution $u$ of problem
\eqref{e1}, which satisfies \eqref{eq4}.
\end{theorem}

\begin{remark}\label{oss1}
Note that in Theorem \ref{thm2} we could replace condition \eqref{eq32} by
\begin{equation}\label{eq32c}
0<\l \leq \frac{(n-1)^2}{4} \underline C^2,
\end{equation}
with $\underline C$ given by Lemma \ref{lemma1}. In this case, the same conclusion is true, replacing $k$ by $\underline C$ in \eqref{eq32} and in \eqref{eq46}. However, $k>0$ is known by assumption, whereas $\underline C$ is not explicitly known (see Remark \ref{oss2}\,).
\end{remark}

Furthermore, if $\g>0$, then for certain values of $\l$ we can improve the assumptions made on $u_0$ in Theorem \ref{thm2}. In fact, we can allow a weaker decaying condition at infinity. In order to do this, we show that, for every
\begin{equation}\label{eq32b}
0<\l < \frac{(n-1)^2}{4} \underline C^2\,,
\end{equation}
with $\underline C>0$ given by Lemma \ref{lemma1}, the function $v$ defined in \eqref{eq16} is a weak supersolution of equation \eqref{e22}
for some $\beta>0$ and for some $\a$ fulfilling
\begin{equation}\label{eq30b}
\max\left\{1-\frac{\g}2, 0\right\}<\alpha<1\,
\end{equation}
(see Proposition \ref{prop3} below). Thus, using this supersolution $v$, we obtain the following result.

\begin{theorem}\label{thm3}
Let assumptions \eqref{A0}-\eqref{A1} be satisfied with $\g>0$. Suppose that
conditions \eqref{eq32b}, \eqref{e25} are satisfied, and that \eqref{eq3} and \eqref{eq2} hold with $w$ replaced by $v$.
Then there exists a global solution $u$ of problem
\eqref{e1}, which satisfies \eqref{eq4}.
\end{theorem}

Finally, in the special case that $\g>2,$ we can further enlarge the class of initial conditions $u_0$. In fact, let $\a>0$. Then we can find $\l^*>0$ such that for every $\l\in (0, \l^*]$, we construct a supersolution $\zeta$ of equation \eqref{e22} such that
\begin{equation}\label{eq34}
\zeta(x)= [r(x)]^{-\alpha}\quad \textrm{for any} \;\; x\in M\setminus B_{R_0}\,,
\end{equation}
for $R_0>0$ big enough (see Proposition \ref{prop4} below). Due to such a supersolution $\zeta$, we get the following result.

\begin{theorem}\label{thm4}
Let assumptions \eqref{A0}-\eqref{A1} be satisfied with $\g>2$. Let $\a>0, \l\in (0, \l^*]$. Suppose that
condition \eqref{e25} are satisfied, and that \eqref{eq3} and \eqref{eq2} hold with $w$ replaced by $\zeta$.
Then there exists a global solution $u$ of problem
\eqref{e1}, which satisfies \eqref{eq4}.
\end{theorem}

Note that in general the solution of problem \eqref{e1} is not unique in $L^{\infty}(M\times (0, T))$. Conditions, related to $M$, that guarantee uniqueness for problem \eqref{e1} are established in \cite{Pu12} (see also e.g. \cite{Grig} for linear equations).

\subsection{Examples}
Set
\[H(t):= \int_0^t h(s) ds\quad \textrm{for any}\;\, t\ge 0\,,\]
and assume that $(A_0)-(A_1)$ are satisfied.

For further references, we recall that in \cite{Pu12} it is proved that if $u_0\ge
0, u_0\not\equiv 0$, and
\begin{equation}\label{e21}
\lim_{t\to \infty} \frac{[H(t)]^{\frac
1{p-1}}}{e^{[\l_1(M)+\e]t}}\,=\,\infty
\end{equation}
for some $\e\in \big(0, \l_1(M)\big)$\,, then any solution to
problem \eqref{e1} blows-up in finite time.

\begin{example}\label{es1}{\em
Assume that $h\equiv 1$ or \eqref{e56} holds. Then we have what follows.

\begin{itemize}
\item{} Let $\g>0$. For any $\l>0$, assumption \eqref{e25} is satisfied. Let $v$ defined by \eqref{eq16} and \eqref{eq30b}. If $u_0\leq C v$, for a sufficiently small $C$, then by Theorem \ref{thm3} we get global existence for problem \eqref{e1}.

\item{} Let $\g\geq 0$. We can apply Theorem \ref{thm2}, but in that case the necessary assumption on the initial datum is worse, since we have to use $v$ with $\a=1$.

\item{} Let $\g>0$. We can apply Theorem \ref{thm1}, if we require that $u_0\leq C w$, for a properly chosen $C>0$. However, concerning the function $w$ we do not have a bound from below, hence the hypothesis on $u_0$ is more implicit in character.

\item{} If $\g>2$, we can further enlarge the class of initial data. In fact, we can suppose that $u_0(x)\leq C [r(x)]^{-\a}$ for any $x\in M\setminus B_{R_0}$, for any $\a>0$, for some  $R_0>0$ and for $C>0$ small enough. Then we apply Theorem \ref{thm4}.
\end{itemize}
}
\end{example}

\begin{example}\label{es2}{\em Let $\s>0$ and
$$h(t)=e^{\s t}\quad\textrm{for any}\;\; t\ge 0.$$
\begin{itemize}
\item{} Suppose that
\begin{equation}\label{eq63}
p > 1 + \frac{\s}{\l_1}\,.
\end{equation}
Then \eqref{e25} holds with $\l=\l_1$. Hence, whenever $\g>0$, we can apply Theorem \ref{thm1}.

\item{} If we replace \eqref{eq63} by a stronger assumption, we can impose conditions on $u_0$ more explicitly. To be specific,
if \begin{equation}\label{eq65}
p > 1 + \frac{\s}{\l}\,,
\end{equation}
for some $\l>0$ satisfying  \eqref{eq32}, then \eqref{e25} is fulfilled for such a $\l$. So, whenever $\g\geq 0$, we can apply Theorem \ref{thm2}, supposing that $u_0\leq C v$ for a suitable $C>0$ with $v$ defined by \eqref{eq16}, \eqref{eq46}, \eqref{eq46b}.

\item{} Furthermore, if \eqref{eq65} is fulfilled for some $\l$ satisfying \eqref{eq32b}, then we can apply Theorem \ref{thm3}, whenever $\g>0$. In this case we must require that $u_0\leq C v$ for a suitable $C>0$ with $v$ defined by \eqref{eq16} and \eqref{eq30b}; so, we can allow a slower decay at infinity on $u_0.$

\item{} In addition, if $\g>2$ we can apply Theorem \ref{thm4}, provided that  \eqref{eq65} is satisfied for $\l>0$ small enough. In this case, we have to impose that $u_0(x)\leq C [r(x)]^{-\a}$ for any $x\in M\setminus B_{R_0}$, for any $\a>0$, for some $R_0>0$ and for $C>0$ small enough.

\item{} Finally, observe that if
\begin{equation}\label{eq64}
p < 1 + \frac{\s}{\l_1}\,,
\end{equation}
then hypothesis
\eqref{e21} is satisfied for appropriate $\e\in \big(0, \l_1(M)\big)$. Hence, we have finite time blow-up.

\item{} We do not know what happens in the case $p=1 + \frac{\s}{\l_1}.$

\end{itemize}}
\end{example}

\begin{example}\label{es3}\em{
We make some comments concerning the relation between our results and those in the literature.
\begin{itemize}
\item{} Let $M=\mathbb H^n$. Hence \eqref{A0}-(ii) holds with $k=1$ and with the equality sign; in addition, $\l_1=\frac{(n-1)^2}{4}$ (see e.g. \cite{Grig4})\,. We can apply Theorem \ref{thm2}; this is in accordance with results in \cite{BPT}. Furthermore, the requested hypothesis on $u_0$ is the same as in \cite{BPT}.

\item{} Assume \eqref{A0} and \eqref{A1}. Hence in order to apply the global existence results in \cite{Pu12}, we must assume that equation \eqref{e22} admits a bounded solution for $\l=\l_1.$ Instead, as we have seen in Examples \ref{es1}-\ref{es2}, we have various global existence results, without making this assumption.

\item{} Assume \eqref{A0} and \eqref{A1}. For $h\equiv 1$, from the results in \cite{Pu14} we get global existence, if $\|u_0\|_{L^{\frac n 2 (p-1)}(M)}$ or $\|u_0\|_{L^p(M)}$  is small enough. This assumption is clearly different in character form those made in Theorems \ref{thm1}, \ref{thm2}, \ref{thm3}, \ref{thm4}. Moreover, it is independent of $\g$, it is indeed only related to assumption \eqref{A1}-(i),(ii).  Moreover, the initial data $u_0$ permitted in Theorems \ref{thm1}, \ref{thm2}, \ref{thm3}, \ref{thm4} not necessarily belong to some $L^p(M)$ space. If $h(t)=e^{\s t}\,\, (\s>0)$, then the the results in \cite{Pu14} require that $p>1 + \frac{\s}{l \l_!},$ for a certain $l=l(p, n)<1$. Hence, this request is worse than that made in  Theorem \ref{thm1}. 

\end{itemize}}
\end{example}

\section{Construction of stationary supersolutions}\label{supersolutions}
In this section we construct supsersolutions to equation \eqref{e22} mentioned in Section \ref{MR}.
First we exhibit the supersolution $w$ for which we do not know the precise behavior at infinity. Then we construct the supersolutions $v$ and $\zeta$.

\smallskip

In the sequel, $v$ is the function defined in \eqref{eq16}.

\subsection{Supersolutions decaying at infinity for $0<\l\leq \l_1$}

\begin{lemma}\label{lemma3}
Let assumption \eqref{A0} be satisfied with $\g>0$. Let $\l>0, \beta>0$ and \eqref{eq30} be satisfied. Then $v$ is a supersolution of equation
\begin{equation}\label{eq17}
\Delta u + \l u = 0\quad \textrm{in} \;\; M\setminus \bar B_{R_0},
\end{equation}
for $R_0>0$ sufficiently large.
\end{lemma}
\begin{proof}
For any $r>0$ we have
\begin{equation}\label{eq18}
v'(r)=-\a \b r^{\a-1} e^{-\b r^{\a}}\,,
\end{equation}
\begin{equation}\label{eq19}
v''(r)=-\a \b e^{-\b r^{\a}}\left[(\a-1) r^{\a-2} - \a \b r^{2\a -2}\right]\,.
\end{equation}
Since $v'<0$, in view of \eqref{e6},  \eqref{eq15}, \eqref{eq30}, \eqref{eq18}, \eqref{eq19}  we deduce that
\begin{equation}\label{eq19b}
\begin{aligned}
& \Delta v(x) + \l v(x) = v''(r) + \mathcal F(r, \theta) v'(r) + \l v(r)  \\
\leq & v''(r) +  \frac{\underline C (n-1)}{r}(1+r)^{1+\frac{\g}{2}} v'(r) +\l v(r) \\
\leq & e^{-\b r^{\a}}\left[\a(1-\a)\b r^{\a-2} + \a^2 \b^2 r^{2\a-2} -\a\b \underline C (n-1) r^{\a-2}(1+r)^{1+\frac{\g}2}+ \l \right]\\
\leq & 0 \quad \textrm{for any}\;\; r\geq R_0\,,
\end{aligned}
\end{equation}
provided that $R_0>0$ is sufficiently large.
\end{proof}

\begin{proposition}\label{prop1}
Let assumption \eqref{A0} be satisfied with $\g>0$. Let $0<\l\leq \l_1, \b>0$. Suppose that \eqref{eq30} holds. Then there exists a weak supersolution $w\in C(M)\cap W^{1,2}_{\textrm{loc}}(M)$ of equation \eqref{e22} such that
\eqref{eq1} is satisfied.
\end{proposition}
\begin{proof}
Let $\phi\in C^{\infty}(M)$ be a positive solution of equation \eqref{e22}; the existence of $\phi$ is guaranteed by results in \cite{CYau} (see also \cite{Grig4}). Let $v$ be the function defined in \eqref{eq16}.
Take $R_1<R_0<R_2,$ with $R_0>0$ given by Lemma \ref{lemma3}. We can find $C>0$ such that
\begin{equation}\label{eq20}
C \phi \leq v \quad \textrm{in} \;\; B_{R_2}\setminus B_{R_1}\,.
\end{equation}
In view of Lemma \ref{lemma3}, by standard results, the function
$$\eta:=\min\{C\phi, v\}$$
is a weak supersolution of equation
\[ \Delta u + \l u = 0 \quad \textrm{in}\;\; M\setminus \bar B_{R_0}\,,\]
in the sense of Definition \ref{solell}.
Define
\begin{equation*}
w(x):=\begin{cases}
\eta &\textrm{in } M\setminus B_{R_0},\\
C\phi &\textrm{in } B_{R_0}.
\end{cases}
\end{equation*}
Due to \eqref{eq20}, it is immediate to see that $w\in C(M)\cap W_{\textrm{loc}}^{1,2}(M)$ is indeed a weak supersolution of equation \eqref{e22}. Finally, from the very definition of $w$ it follows that \eqref{eq1} holds.
\end{proof}

\subsection{Further supersolutions for small $\l>0$}
\begin{proposition}\label{prop2}
Let assumption \eqref{A0} and conditions \eqref{eq32}, \eqref{eq46}, \eqref{eq46b} be satisfied.  Then $v\in C(M)\cap W^{1,2}_{\textrm{loc}}(M)$ is a weak supersolution of equation \eqref{e22}.
\end{proposition}

\begin{proof} For every $r>0$ we have
\begin{equation}\label{eq21}
v'(r)=-\b e^{- \b r}\,,
\end{equation}
\begin{equation}\label{eq22}
v''(r)=\b^2 e^{- \b r}\,.
\end{equation}
Since $v'<0$, in view of \eqref{e6}, \eqref{eq15b}, \eqref{eq21}, \eqref{eq22}  we deduce that
\begin{equation}\label{eq23}
\begin{aligned}
 \Delta v(x) + \l v(x)  = & v''(r) + \mathcal F(r, \theta) v'(r) + \l v(r)  \\
\leq & v''(r) + k v'(r) +\l v(r) \\
\leq & e^{-\b r}\left[\b^2 - \b k (n-1) +\l \right]\leq  0 \quad \textrm{in}\;\; M\setminus \{o\}\,.
\end{aligned}
\end{equation}
Here we have used hypotheses \eqref{eq32} and \eqref{eq46b}.
Now, recall that the function $x\mapsto r(x)$ is $C^2$ in $M\setminus\{o\}$ and $1-$Lipschitz in the whole of $M$. Thus, we can infer that in the weak sense
\begin{equation}\label{eq24}
\Delta v + \l v \leq 0 \quad \textrm{in}\;\; M\,.
\end{equation}
\end{proof}

\begin{proposition}\label{prop3}
Let assumption \eqref{A0} be satisfied with $\g>0$. Assume that condition \eqref{eq32b} holds. Then, for some $\b>0$ and for some $\a$ fulfilling \eqref{eq30b}, $v\in C(M)\cap W^{1,2}_{\textrm{loc}}(M)$  is a bounded weak supersolution of equation \eqref{e22}.
\end{proposition}

\begin{proof} Consider the function $v$ defined in \eqref{eq16}. From \eqref{eq19b}, we obtain that, for all $x\equiv (r, \theta)\in M\setminus \{o\},$
\begin{equation}\label{eq36}
\Delta v(x) +\l v(x) \leq e^{-\b r^{\a}} r^{\a-2} \varphi(r)\,,
\end{equation}
where
\begin{equation}\label{eq70}
\varphi(r):= \a(1-\a)\b + \a^2 \b^2 r^{\a} -\a\b \underline C (n-1)(1+r)^{1+\frac{\g}2}+ \l r^{2-\a}\,.
\end{equation}
Due to \eqref{eq30b}, for any $r\in (0,1]$ we obtain
\begin{equation}\label{eq37}
\varphi(r) \leq \a^2\b^2 + [\a-\a^2 -\a\underline C(n-1)] \b + \l\,.
\end{equation}
By hypothesis \eqref{eq32b}, we can find $\a$ such that \eqref{eq30b} holds and
\begin{equation}\label{eq38}
\l\leq \frac 1 4[1-\a- \underline C(n-1)]^2\,.
\end{equation}
Hence, for some $\b>0$,
\begin{equation}\label{eq39}
\a^2\b^2 + [\a-\a^2 -\a\underline C(n-1)] \b + \l \leq 0\,.
\end{equation}
Thus, in view of \eqref{eq37} and \eqref{eq39} we get
\begin{equation}\label{eq40}
\varphi(r) \leq 0 \quad \textrm{for every}\;\; r\in (0, 1]\,.
\end{equation}
Observe that, for this choice of $\a$ and $\b$, we have
\begin{equation}\label{eq41}
\a^2\b^2 -\a\b \underline C(n-1) +\l \leq -\a(1-\a)\b<0\,.
\end{equation}
Furthermore, for any $r>1$, due to \eqref{eq30b}, \eqref{eq39} and \eqref{eq41} we have
\begin{equation}\label{eq41b}
\begin{aligned}
\varphi(r)\leq & \a(1-\a)\b + r^{2-\a}[\a^2 \b^2 -\a\b \underline C(n-1) +\l] \\ \leq &  \a^2\b^2 + [\a-\a^2 -\a\underline C(n-1)] \b + \l\ \leq 0\,.
\end{aligned}
\end{equation}
From \eqref{eq36}, \eqref{eq40} and \eqref{eq41b} we can deduce that
\begin{equation}\label{eq42}
\Delta v + \l v \leq 0   \quad \textrm{in}\;\; M\setminus \{o\}\,.
\end{equation}
Indeed, $v$ is a weak supersolution of equation \eqref{e22}, by the same arguments used in the proof of Proposition \ref{prop1}.
\end{proof}

We have also the following result.
\begin{proposition}\label{prop5}
Let assumption \eqref{A0} be satisfied. Assume that condition \eqref{eq32c} holds, and that
\begin{equation}\label{eq71}
1\leq \a \leq \min\left\{1+\frac{\g}2 , 2 \right\}\,.
\end{equation}
Then, for some $\b>0$, $v\in C(M)\cap W^{1,2}_{\textrm{loc}}(M)$  is a bounded weak supersolution of equation \eqref{e22}.
\end{proposition}
\begin{proof}
We have that \eqref{eq36} holds with $\varphi$ defined in \eqref{eq70}. In view of \eqref{eq71}, for any $r\in (0, 1]$ we obtain
\begin{equation}\label{eq37c}
\varphi(r) \leq \a^2\b^2 - \a\underline C(n-1) \b + \l\,.
\end{equation}
Moreover, due to \eqref{eq71} again, we can infer that for any $r>1$
\begin{equation}\label{eq72}
\varphi(r) \leq r^{1+\frac{\g}2}[\a^2\b^2 - \a\underline C(n-1) \b + \l]\,.
\end{equation}
By hypothesis \eqref{eq32c}, we can find $\b>0$ such that
\[\a^2\b^2 - \a\underline C(n-1) \b + \l\leq 0.\]
So, from \eqref{eq37c} and \eqref{eq72} we have that
\[\varphi(r)\leq 0 \quad \textrm{for every}\;\; r>0\,.\]
Hence the conclusion follows as in the proof of Proposition \ref{prop4}.
\end{proof}

\begin{remark}\label{oss4} Proposition \ref{prop5} is given for the sake of completeness. However, in the global existence results, we do not consider the supersolution $v$ given by Proposition \ref{prop5}, since that given by Proposition \ref{prop3} allows us to deal with a larger class of initial conditions, provided \eqref{eq32b} holds. On the other hand, let $\l=\frac{(n-1)^2}4 \underline C^2$. The best choice of $\a$ in Proposition \ref{prop5}, with respect to the class of permitted initial data $u_0$, is $\a=1$, which has been already considered in Remark \ref{oss1}.
\end{remark}

\begin{proposition}\label{prop4}
Let assumption \eqref{A0} be satisfied with $\g>2$. Let $\a>0$. Then, for some $\l^*>0$, for every $0<\l\leq \l^*$ there exists a weak supersolution $\zeta\in C(M)\cap W^{1,2}_{\textrm{loc}}(M)$ of equation \eqref{e22} such that \eqref{eq34} is satisfied.
\end{proposition}
\begin{proof}
Define
\[\zeta_1(x)\equiv \zeta_1(r):=[r(x)]^{-\alpha}\quad \textrm{for all} \;\; x\in M\setminus \{o\}\,.\]
For any $r>0$ we have
\[\zeta_1'(r)=-\a r^{-\a-1}\,,\]
\[\zeta_1''(r)=\a(\a+1)r^{-\a-2}\,.\]
Thus, since $\zeta_1'<0$, in view of \eqref{e6}, \eqref{eq15}, \eqref{eq30}, using condition $\g>2$, we get
\begin{equation}\label{eq50}
\begin{aligned}
& \Delta \zeta_2(x) + \l \zeta_2(x)  \\  & \leq  \a(\a+1) r^{-\a-2} -\a \underline C(n-1) r^{-\a-1+\frac{\g}{2}} + \l r^{-\a} \\
& =  r^{-\a}[\a(\a+1)r^{-2} -\a\underline C (n-1) r^{-1+\frac{\g}2}+\l]\leq 0\quad \textrm{in}\;\, M\setminus B_{R_0},
\end{aligned}
\end{equation}
provided that $R_0>0$ is large enough.

\smallskip

Now, for each $a>0, b>0,$ define
\[\zeta_2(x)\equiv \zeta_2(r):=a - b r(x)  \quad \textrm{for all } \; x\in M\,.\]
We choose $a>0$ and $b>0$ so that
\begin{equation}\label{eq51}
\zeta_1(R_0)=\zeta_2(R_0)\,, \;\; \zeta_2'(R_0)\geq \zeta_1'(R_0)\,.
\end{equation}
It is easily checked that if we take
\[b= \a R_0^{-\a-1}\,,\;\; a=b R_0+ R_0^{-\a}\,,\]
then conditions in \eqref{eq51} hold; in addition,
\[\zeta_2 > 0 \quad \textrm{in}\;\; \bar B_{R_0}\,.\]
Furthermore, in view of \eqref{eq15b}, we have
\[
\begin{aligned}
\Delta \zeta_2(x) + \l \zeta_2(x)& \\ & \leq - b k(n-1) + \l a -\l b r \\ &\leq -\a R_0^{-\a-1}k(n-1) + \l (\a+1) R_0^{-\a}\leq 0 \quad \textrm{in} \;\, B_{R_0}\setminus\{o\}\,,
\end{aligned}
\]
provided $\l>0$ is small enough. By the same arguments as in the proof of Proposition \ref{prop1}, indeed $\zeta_2$ solves weakly
\begin{equation}\label{eq52}
\Delta \zeta_2(x) +\l \zeta_2(x) \leq 0 \quad \textrm{in}\;\; B_{R_0}\,.
\end{equation}

Define
\[\zeta:=\begin{cases}
\zeta_1 &\textrm{in } M\setminus B_{R_0},\\
\zeta_2 &\textrm{in } B_{R_0}.
\end{cases}
\]
From \eqref{eq50}, \eqref{eq51} and \eqref{eq52}, by standard tools, it is immediately seen that $\zeta$ is a weak supersolution of equation
\eqref{e22}. Clearly, \eqref{eq34} holds.
\end{proof}

\section{Global existence: proofs}\label{dim1} \setcounter{equation}{0}
\begin{proof}[Proof of Theorem \ref{thm1}\,.] Let $\{\Omega_j\}_{j\in\ene}$ be a sequence of domains
$\{\Omega_j\}_{j\in \ene}\subseteq M$ such that
$\bar\Omega_j\subseteq \Omega_{j+1}$ for every $j\in \ene,\,
\bigcup_{j=1}^\infty\Omega_j=M\,,\pa \Omega_j$ is smooth for every
$j\in \ene\,.$

\smallskip

For any $j\in \ene$ let $u_j$ be the unique classical solution to
problem
\begin{equation}\label{e28}
\left\{
\begin{array}{ll}
 \,  \pa_t u = \Delta u \,+ h(t) u^p &\textrm{in}\,\,\Omega_j\times (0,T)
\\&\\
\textrm{ }u \,= 0 & \textrm{in}\,\, \pa\Omega_j\times (0,T) \\&\\
\textrm{ }u \, = u_0& \textrm{in}\,\,  \Omega_j\times \{0\} \,.
\end{array}
\right.
\end{equation}

\smallskip
Take the constant $\tilde C>0$ given by \eqref{eq2}. Let
\[\tilde w:= \tilde C w\quad \big(x\in M\big)\,,\]
and
\[\xi(t)=[1-(p-1)\|\tilde w \|_{\infty}^{p-1}\tilde H(t)]^{-\frac 1{p-1}}\,\quad \big(t\in [0,\infty)\big)\,;\]
note that $\xi$ is well-defined in $[0, \infty)$ due to \eqref{e25} and \eqref{eq2}.
It is easily seen that $\xi$ solves
problem
\begin{equation}\label{e39}
\left\{
\begin{array}{ll}
 \,  \xi' = \|\tilde w \|_\infty^{p-1}\tilde h(t)\xi^p\,,  &
 t\in (0,\infty)
\\&\\
\textrm{ }\xi(0) \, = 1 \,.
\end{array}
\right.
\end{equation}
Define
\[\bar u(x,t):= e^{-\l t}\xi(t)\tilde w (x)\quad \big((x,t)\in M\times[0,\infty)\big)\,. \]
In view of Proposition \ref{prop1} and \eqref{e39}, for each $0<\l\leq \l_1$, we have 
\[\partial_t \bar u - \Delta \bar u - h(t) \bar u^p \]\[= -\l e^{-\l t}\xi(t) \tilde w(x) + e^{-\l t}\|\tilde w \|_{\infty}^{p-1}h(t) e^{-(p-1)\l t} \xi^p(t) \tilde w(x) \]\[+ \l e^{-\l t}\xi(t) \tilde w(x) -h(t) e^{-\l p t}\xi^p(t) \tilde w^p (x) \geq 0\quad  \textrm{weakly in}\,\, M\times (0, \infty)\,. \]
So, $\bar u$ is a weak supersolution of equation
\begin{equation}\label{eq67}
\pa_t u = \Delta  u + h(t)  u^p \quad \textrm{in} \;\; M\times (0, \infty)\,.
\end{equation}
Moreover, due to \eqref{eq3}, we have that for any $j\in \ene$, $\bar u$ is a bounded weak
supersolution of problem \eqref{e28}. Obviously, for each $j\in \ene$, $\underline u\equiv 0$ is a subsolution of problem \eqref{e28}.
Hence, by the comparison principle for every $j\in \ene$ we obtain
\begin{equation}\label{e29}
0 \le u_j \le \bar u\quad \textrm{in}\;\; B_j\times (0, T)\,.
\end{equation}
In view of \eqref{e29}, by standard compactness arguments (see e.g. \cite{Fr}, \cite{LSU}, \cite{L}), there exists a subsequence
$\{u_{j_k}\}\subseteq \{u_j\}$, which converges locally uniformly
in $M\times (0,T)$ to a weak solution $u$ of problem \eqref{e1}. By standard regularity results (see e.g. \cite{Fr}, \cite{LSU}, \cite{L}) , indeed $u$ is a classical solution of equation \eqref{eq67}; moreover, by using local barrier arguments (see e.g. \cite{Fr}) it follows that  $u\in C(M\times [0, T))$ and $u=u_0$ in $M\times \{0\}$.
Furthermore, from \eqref{e29} we get
\[0\le u\le \bar u\quad\textrm{in}\;\;M\times (0, T)\,.
\]
Hence the thesis follows.
\end{proof}

Clearly, if $h$ and $u_0$ are Holder continuous, then, in view of standard methods (see e.g. \cite{Fr}, \cite{LSU}, \cite{L}), we have that $u_{j_k}$ converges locally in $C^{2,1}_{x,t}(M\times [0, T)$ to a solution of problem \eqref{e1}.

\smallskip

Theorems \ref{thm2}, \ref{thm3} and \ref{thm4} can be proved by the same arguments as in the proof of Theorem \ref{thm1}, replacing
$w$ by $v$ in the proofs of Theorems \ref{thm2}, \ref{thm3}, and by $\zeta$ in the proof of Theorem \ref{thm4}, and choosing $\l$ appropriately. 

\begin{remark}\label{oss3} In \cite{Pu12} in the proof of global existence it is used a comparison principle in $M\times (0, T)$; this holds under suitable assumptions on $M$.
Indeed, here we use comparison principles only in $B_j\times (0, T)$ for every $j\in \ene$, hence we do not require those assumptions. Moreover, observe that in particular in Theorem \ref{thm4} when \eqref{A0} is satisfied with $\g>2$, the comparison principle on $M\times (0, T)$ does not hold.
\end{remark}

\end{document}